\documentclass[reqno,11pt]{amsart}
\usepackage{amsmath,amsthm,amssymb}
\usepackage{graphicx}
\usepackage{changes}

\usepackage{xcolor}

\usepackage[section]{placeins}

\def\diff{\mathrm{d}}

\newtheorem{lemma}{Lemma}[section]
\newtheorem{theorem}[lemma]{Theorem}

\newtheorem{proposition}[lemma]{Proposition}
\theoremstyle{remark}

\numberwithin{equation}{section}

\title[Threshold for blowup for the  cubic wave equation]{Threshold for blowup for the  supercritical cubic wave equation}

\author{Irfan Glogi\'c}
\address{Universit\"at Wien, Fakult\"at f\"ur Mathematik,
	Oskar-Morgenstern-Platz 1, A-1090 Vienna, Austria}
\email{irfan.glogic@univie.ac.at}

\author{Maciej Maliborski}
\address{Gravitational Physics, Faculty of Physics, University of Vienna, Boltzmanngasse 5, A-1090 Vienna, Austria} \email{maciej.maliborski@univie.ac.at}

\author{Birgit Sch\"orkhuber}
\address{Karlsruhe Institute of Technology, Institute for Analysis,  Englerstra{\ss}e 2, 76131 Karlsruhe, Germany} \email{birgit.schoerkhuber@kit.edu}

\thanks{Irfan Glogi\'c is supported by the Austrian Science Fund FWF,
	Project P 30076. Birgit Sch\"orkhuber gratefully acknowledges support from the Klaus-Tschira Stiftung.
	Computations have been performed on Minerva cluster of the
	Max-Planck Institute for Gravitational Physics. Funded by the Deutsche Forschungsgemeinschaft (DFG, German Research Foundation) - Project-ID 258734477 - SFB 1173}

\begin{document}
\begin{abstract}
We consider the focusing cubic wave equation in the energy supercritical case, i.e., in dimensions $d \geq 5$. For this model an explicit nontrivial self-similar blowup solution was recently found by the first and third author in \cite{GlogicSchoerkhuber}. Furthermore, the solution is proven to be co-dimension one stable in $d=7$.  In this paper, we study the equation from a numerical point of view. For $d=5$ and $d=7$ in the radial case, we provide evidence that this solution is at the threshold between generic ODE blowup and dispersion. In addition, we investigate the spectral problem that underlies the stability analysis  and compute the spectrum in general supercritical dimensions.
\end{abstract}

\maketitle

\section{Introduction}
\label{sec:1}

We consider the focusing cubic wave equation
\begin{equation}\label{Eq:NLW}
( \partial^2_t -  \Delta )u(t,x)   =   u(t,x)^{3}   \\
\end{equation}
$(t,x) \in I \times \mathbb R^d$, where $I \subset \mathbb R$ is an interval containing zero. The model is invariant under rescaling $u \mapsto u_{\lambda}$
\begin{equation*}
u_{\lambda}(t,x) = \lambda^{-1} u(t/\lambda , x/\lambda ), \quad \lambda > 0.
\end{equation*}
The corresponding scale invariant Sobolev space for $(u(t,\cdot), \partial_t u(t,\cdot))$ is  $\dot H^{s_c} \times \dot H^{s_c-1}(\mathbb R^d)$, $s_c = \frac{d}{2}-1$.
Note that  $s_c= 1$ for $d=4$ and in this case, Eq.~\eqref{Eq:NLW} is referred to as energy critical.  In this paper, we restrict ourselves to the supercritical case $d \geq 5$. Eq.~\eqref{Eq:NLW} admits a trivial self-similar blowup solution,
\begin{equation}\label{Eq:ODEblowup}
u_T(t,x)  =  \frac{\sqrt{2}}{T-t}, \quad T > 0,
\end{equation}
which is known to be stable, see  \cite{ChatzDonn2019}, \cite{DonningerSchoerkhuber2017}. In the recent work \cite{GlogicSchoerkhuber} by the first and the third author a non-trivial self-similar solution was found
\begin{equation}\label{Eq:CritSol}
u^*_T(t,x)=\frac{1}{T-t} U^*\left(\tfrac{|x|}{T-t}\right), \quad U^*(\rho) =  \frac{2 \sqrt{2(d-1)(d-4)}}{d-4+3 \rho^2}
\end{equation}
for $d \geq 5$.  Furthermore, in $d=7$ this solution is proven to be co-dimension one stable. More precisely, there exists a co-dimension one Lipschitz manifold  of initial data in a small neighbourhood of $u^*_T$, whose solutions blow up in finite time and converge asymptotically to $u^*_T$ (modulo space-time shifts and Lorentz boosts) in the backward light cone of the blowup point.

In Sec.~\ref{Codimension_1_stability},  we address the stability of $u^*_T$  and numerically investigate the underlying spectral problem for general dimensions.
In Sec.~\ref{sec:Non-genericBlowup}, we provide evidence that this solution is an attractor within a co-dimension one manifold that is a threshold between ODE blowup and dispersion. To finish the introduction, we give a short overview on threshold phenomena in energy supercritical wave equations, which is our main motivation to study the model at hand.

\subsection{Threshold for blowup in energy supercritical models} \label{Sec:Threshold}

In the past decades, a vast body of literature has been concerned with the question of global existence versus finite-time blowup of solutions in nonlinear PDEs. 
Intuitively, one expects  ``small" initial data to lead to global in time solutions, while singularities are supposed to form from sufficiently ``large" data. Consequently, one could expect the existence of a certain threshold between these two basins of attraction. From a mathematical point of view, this is highly non-trivial to be made precise; in particular in energy supercritical models, where little is known in general.

In the context of gravitational collapse, this question has been investigated numerically in the 90s  by Choptuik  \cite{Choptuik1993}  for a simple matter model. By considering one-parameter families of solutions interpolating between dispersion and black hole formation, he found a (discretely) self-similar solution as an intermediate attractor for solutions close to the threshold between these two end states. From a physical point of view, this solution violates conjectures about the structure of gravitational singularities. In the past years, numerous simulations have investigated so-called critical phenomena in more involved models, see e.g.~\cite{GunMar07} for a review. It was found that threshold solutions between different stable regimes are either static or self-similar (discretely or continuously). Remarkably, the existence of the Choptuik critical solution has recently been proven by Reiterer and Trubowitz in \cite{MR3946407} by using computer-assisted methods.

In the recent years it has become evident that self-similar threshold solutions seem to be a feature of energy supercritical models (including Einstein's equation) rather than being specific to gravitational collapse: for the wave maps equation,  a critical self-similar blowup solution was observed numerically in $d=3$ by Bizo\'n, Chmaj and Tabor \cite{BizonChmajTabor2000} and by Biernat,  Bizo\'n and the second author in higher space dimensions \cite{BiernatBizonMaliborski2016}. Other examples are the Yang-Mills equation on $\mathbb R^{5+1}$ \cite{Bizon2002_YM}, and the supercritical focusing wave equation in three dimensions \cite{BizonChmajTabor2004}. 
From an analytic point of view, these problems are widely open. In fact, the model considered in this paper seems to be the first wave equation for which a candidate for a critical self-similar solution is known in closed form. Thus, it provides a good starting point for the analytic investigation of threshold phenomena in supercritical wave equations. 

We remark in passing that critical self-similar solutions have been observed numerically in other problems such as the three-dimensional parabolic-elliptic Keller-Segel model \cite{MR1709861}. There, the situation seems to be even more complex due to the existence of different stable blowup regimes and at least two different critical self-similar solutions.

\subsection{Blowup dynamics for the supercritical cubic wave equation}\label{Sec:KnownResults}

To complete our discussion, we mention some known results about singularity formation for the model under investigation. For a more thorough discussion on the focusing wave equation, we refer the reader e.g.~to \cite{DonningerSchoerkhuber2017}.

For Eq.~\eqref{Eq:NLW}, the generic blowup behavior is conjectured to be governed by the ODE blowup solution \eqref{Eq:ODEblowup}. In fact, this is also corroborated by  our numerical simulations presented in Sec.~\ref{sec:Non-genericBlowup}. From a rigorous point of view, the stability of this solution (in a local sense) has been proved by Donninger and the third author in \cite{DonningerSchoerkhuber2017} for all odd $d \geq 5$ in spherical symmetry. The non-radial case was addressed by Donninger and Chatzikaleas \cite{ChatzDonn2019} in $d=5$.

Non-trivial self-similar solutions have been investigated numerically by Kycia \cite{MR2804023}; from that one can expect the existence of infinitely many profiles $\{U_n : n \in \mathbb N_0 \}$  in dimensions $5 \leq d < 13$, with $U_0 = \sqrt{2}$.  In  $d \geq 13$, a result by Collot \cite{Collot2018} proves the existence of non-self-similar blowup solutions with more than two unstable directions.

We conclude this section by briefly commenting on the  energy critical case. There, the threshold is characterized in terms of the static ground state solution and this is fairly well understood from an analytic point of view. Results for Eq.~\eqref{Eq:NLW} in $d=4$ have been obtained for example in \cite{MR3769743},  \cite{DuyckaertsMerle2008}, \cite{MR3006642} and  \cite{KenigMerle2008}. We also refer the reader to \cite{KriegerNakanishiSchlag2015} for a characterization of the  threshold for an energy critical wave equation in three dimensions. However, these techniques are specific to the critical case and cannot be transferred directly to supercritical problems.

\section{Co-dimension one stable self-similar blowup} \label{Codimension_1_stability}


To investigate the stability of $u_T^*$ in the spherically symmetric setting,  we fix $d \geq 5$ and consider the radial cubic wave equation
\begin{equation}
\label{eq:1}
\partial_{t}^{2}u(t,r)- \partial_{r}^{2}u(t,r) - \frac{d-1}{r}\partial_{r}u(t,r) = u(t,r)^{3}
\end{equation}
for small radial perturbations $(f,g)$ of the blowup initial data
\begin{equation}\label{Perturbed_blowup_data}
u(0,\cdot ) = u^*_1(0,\cdot)  + f, \quad \partial_t u(0,\cdot) = \partial_t u^*_1(0,\cdot)  + g.
\end{equation}
Since we are interested in the blowup behavior near the origin, we consider the evolution in backward light cones 
\begin{equation}\label{Lightcone}
\mathcal C_{T} :=\bigcup_{t \in [0,T)} \{t\} \times \overline{\mathbb  B^{d}_{T-t}}.
\end{equation}

We look for solutions that can be written as
\begin{equation}\label{Ansatz_solution}
u(t,r) =  u^*_T(t,r) + (T-t)^{-1} \varphi(- \log(T-t) , \tfrac{r}{T-t} )
\end{equation}
for some $T > 0$, such that the perturbation $\varphi$ vanishes in a suitable sense as $t \to T^{-}$. The result of \cite{GlogicSchoerkhuber} proves this in $d=7$ under a co-dimension one condition on the initial data.

\begin{theorem}[\cite{GlogicSchoerkhuber}, radial version] \label{Codim_one_radial}
	Let $d=7$ and define
	\begin{equation}\label{Eq:Def_UnstableDir}
f_1(r) = (1+r^2)^{-2} , \quad g_1(r) = 4 (1+r^2)^{-3}.
\end{equation}
There are constants $\omega, \delta,c > 0$ such that for all smooth, radial $(f,g)$ with
	\begin{equation*}
	\| (f,g) \|_{H^4 \times H^3(\mathbb B_2^7)} \leq \tfrac{\delta}{c}
	\end{equation*}
	the following holds: There are $\alpha \in[-\delta,\delta]$ and  $T  \in [1-\delta, 1+ \delta]$ depending Lipschitz continuously on $(f,g)$ such that for initial data
	\begin{equation}\label{Codim_one_data}
	u(0,\cdot ) = u^*_1(0,\cdot)  + f + \alpha f_1, \quad \partial_t u(0,\cdot) = \partial_t u^*_1(0,\cdot)  + g + \alpha g_2
	\end{equation}
	there is a unique solution $u$ in the backward light cone $\mathcal C_T$ blowing up at $t = T$ and converging to $u^*_T$ according to
	\begin{eqnarray}\label{MainTh_bounds}
	\begin{split}
	(T-t)^{k-s_c}  \| u(t,\cdot) - u^*_{T}(t,\cdot) \|_{\dot H^{k} (\mathbb B_{T-t}^7)}\lesssim (T-t)^{\omega}  \\
	(T-t)^{k-s_c}  \| \partial_t u(t,\cdot) - \partial_t u^*_{T}(t,\cdot) \|_{\dot H^{k-1} (\mathbb B_{T-t}^7)} \lesssim (T-t)^{\omega}
	\end{split}
	\end{eqnarray}
	for $k=1,2,3$.
\end{theorem}

We note that the right-hand side of Eq.~\eqref{MainTh_bounds} is normalized to the behavior of $u^*_T$ in the respective norm.
Furthermore, by a solution, we mean a solution to the corresponding operator equation in adapted coordinates, see \cite{GlogicSchoerkhuber}.
The co-dimension one condition is formulated in terms of explicit functions $(f_1,g_1)$ which arise as eigenfunctions in the spectral problem for the linearization around the blowup solution. This will be explained in more detail below. 

The proof of Theorem \ref{Codim_one_radial} relies on the analysis of the time evolution for the perturbation $\varphi$ in self-similar coordinates
\begin{equation}\label{SS_coordinates}
\tau = - \log(T-t), \quad \rho = \frac{r}{T-t}.
\end{equation}
Note that the backward light cone is mapped to $\bigcup_{\tau \in [0,\infty)} \{\tau\} \times [0,1]$.
By defining
\begin{equation}\label{SS_coordinates2}
\psi(\tau,\rho) := e^{-\tau} u(T-e^{-\tau}, e^{-\tau} \rho ) = U^*(\rho) + \varphi(\tau,\rho),
\end{equation}
the evolution for $\varphi$ is given by
\begin{equation}\label{Evol_perturb}
\left(\partial^2_{\tau} + 3 \partial_{\tau} + 2 \rho \partial_{\rho} \partial_{\tau} - \Delta_{\rho} + \rho^2 \partial^2_{\rho} + 4 \rho \partial_{\rho} + 2 - V(\rho) \right ) \varphi(\tau,\rho)  = N(\varphi(\tau,\rho))
\end{equation}
with $V(\rho) = 3 U^*(\rho)^2$ and $N(\varphi) = (U^* + \varphi)^3 -  3 {U^*}^2$.  It can be shown that the dynamics are governed by the linearized problem and that the nonlinearity $N$ on the right-hand side can be treated perturbatively.  The key to a generalization of Theorem \ref{Codim_one_radial} to other (odd) space-dimensions  is the solution of the spectral problem for the corresponding linearization.  As a matter of fact,  it suffices to study \textit{mode solutions}
\begin{equation}\label{ModeAnsatz}
\varphi(\tau,\rho) = e^{\lambda \tau} f(\rho)
\end{equation}
with $\lambda \in \mathbb C$, $\mathrm{Re} \lambda \geq 0$ and \textit{smooth} profiles $f \in C^{\infty}[0,1]$, that satisfy the linearized equation
\begin{equation}\label{Evol_perturb_lin}
\left(\partial^2_{\tau} + 3 \partial_{\tau} + 2 \rho \partial_{\rho} \partial_{\tau} - \Delta_{\rho} + \rho^2 \partial^2_{\rho} + 4 \rho \partial_{\rho} + 2 - V(\rho) \right ) \varphi(\tau,\rho)  =0.
\end{equation}

In a rigorous formulation, the values $\lambda$ correspond to eigenvalues of a suitably defined differential operator $\mathbf L_0$, see \cite{GlogicSchoerkhuber}.

\subsection{The spectral problem}\label{Sec:ODE_Analysis}
A fundamental observation is that the time translation invariance of the blowup solution induces a mode solution for $\lambda_0 = 1$. More precisely,
\begin{equation}\label{Gauge_mode_phys}
\phi_0(t,r) := \frac{d}{d \varepsilon} u^*_{T+\varepsilon}(t,r)|_{\varepsilon=0} =(T-t)^{-2} f_0(\tfrac{r}{T-t})
\end{equation}
with
\begin{equation}\label{gauge_mode}
f_0(\rho) =  \frac{d-4-3\rho^{2}}
{\left(d-4+3\rho^{2}\right)^{2}},
\end{equation}
satisfies the linearized equation in physical variables $(t,r)$,
\begin{equation}\label{Lin_Eq_phys}
\left (\partial_{t}^{2} - \partial_{r}^{2} - \frac{d-1}{r}\partial_{r} - 3 (T-t)^{-2}  U^*(\tfrac{r}{T-t})^{2} \right) \phi(t,r) = 0.
\end{equation}
Hence, $\varphi_0(\tau,\rho) = e^{-\tau} \phi_0(T-e^{-\tau}, e^{-\tau} \rho ) =  e^{\tau} f_0(\rho)$ is a solution of  \eqref{Evol_perturb_lin}.
In the nonlinear time-evolution this gauge mode can be controlled by suitably adjusting the blowup time $T >0$ and is not an obstruction to stability. However, as will be discussed in the following, there is an additional genuine instability that yields the co-dimension one condition on the data.


 By inserting the mode ansatz \eqref{ModeAnsatz} into Eq.~\eqref{Evol_perturb_lin} we obtain the following ordinary differential equation
 \begin{multline}\label{Eq:ODE_L0}
 (1-\rho^2) f''(\rho) +  \left(\frac{d-1}{\rho}   - 2(\lambda +2) \rho \right) f'(\rho)    \\
 -\left( (\lambda + 1)(\lambda + 2)  - \frac{24(d-1)(d-4)}{(d-4+3 \rho^2)^2} \right )f(\rho)=0.
 \end{multline}
 We are interested in values of  $\lambda$ that yield solutions that are smooth on $[0,1]$. Note that \eqref{Eq:ODE_L0} is a Fuchsian equation with six singular points. In particular, $\rho=0$ and $\rho=1$ are singularities and Frobenius theory implies that eigenfunctions are not just smooth but analytic on $[0,1]$. This observation allows for an effective use of the shooting method to compute the eigenvalues.
 	This standard technique relies on approximately computing analytic solutions emanating from two singular points and then adjusting the underlying parameter so that solutions smoothly match at some third, conveniently chosen point. Our particular approach follows the work of Biernat and Bizo\'n \cite{BizonBiernat2015} on an analogous problem relative to supercritical wave maps.
  
 	To simplify the analysis we first reduce the number of singular points of Eq.~\eqref{Eq:ODE_L0}. This is done by the change of the independent variable
 	\begin{equation}\label{Change_var}
 	x=\rho^2.
 	\end{equation}
 	In this way, the set of (regular) singular points becomes $\{0,1,(4-d)/3,\infty\}$. We remark that Fuchsian equations with four singular points go under the name of Heun, and to bring the equation to its canonical Heun form, see \cite{NIST10}, we also scale the dependent variable
 	\begin{equation}\label{Change_var_Heun}
 	f(\rho)=\frac{y(x)}{(d-4+3x)^2}.
 	\end{equation}
 	In this way we are lead to the following equation
 \begin{multline}\label{Heun}
 y''(x)+\left(\frac{d}{2x}+\frac{2\lambda+5-d}{2(x-1)}-\frac{12}{3x+d-4} \right)y'(x)\\+\frac{3(\lambda-3)(\lambda-2)x+\lambda(\lambda+3)(d-4)-10d+16)}{4x(x-1)(3x+d-4)}y(x)=0.
 \end{multline}
 Since the set of Frobenius indices of Eq.~\eqref{Eq:ODE_L0} at $\rho=0$ is $\{0,2-d\}$ analytic solutions at $\rho=0$ are even and therefore \eqref{Change_var} preserves analyticity at the origin. What is more, Eq.~\eqref{Eq:ODE_L0} and Eq.~\eqref{Heun} have the same set of eigenvalues.
 
 The sets of Frobenius indices for Eq.~\eqref{Heun} at $x=0$ and $x=1$ are $\{0,1-d/2\}$ and $\{0,(d-3)/2-\lambda\}$ respectively. 
 Therefore, if $(d-3)/2-\lambda$ is not a positive integer then the (normalized) analytic solutions at $x=0$ and $x=1$ are
 	\begin{equation}\label{series}
 	y_0(x)=1+\sum_{n=1}^{\infty}a_nx^n \quad \text{and} \quad
 	y_1(x)=1+\sum_{n=1}^{\infty}b_n(1-x)^n 
 	\end{equation}
 	respectively. Note that the scenario that yields an eigenvalue is precisely when $y_0$ and $y_1$ are constant multiples of each other, which is in turn equivalent to Wronskian
 \begin{equation}\label{wronskian}
 W[y_0,y_1](x):=y_0'(x)y_1(x)-y_0(x)y_1'(x)
 \end{equation}
 being identically zero. Now we approximately compute the solutions $y_0$ and $y_1$ by truncating the series \eqref{series} for a large $n$ and then we evaluate the Wronskian \eqref{wronskian} at the midpoint $x=1/2$. Of course, the choice of this point may depend on other singularities, as their position (in general) determines the radius of convergence of series \eqref{series}.   The (approximate) eigenvalues are then given by the zeroes of the Wronskian. We remark that for the Heun equation there is no closed form expression for the Wronskian, unlike for example for the hypergeometric equation where \eqref{wronskian} is given in terms of Gamma functions. However, functions the $y_0$ and $y_1$ are built into Maple, and we can therefore numerically compute the zeros of $W[y_0,y_1](1/2)$ with relatively high precision.
 
 Still, recall that the shooting is done with the assumption that $(d-3)/2-\lambda$ is not a positive integer. Therefore, this complementary case has to be treated separately. Namely, in this situation it can happen that both Frobenius solutions at $x=1$ are analytic, in which case the underlying $\lambda$ is obviously an eigenvalue. Otherwise only the subdominant Frobenius solution is analytic and by factoring out its asymptotic behavior at $x=1$ we obtain a new equation which is amenable to the shooting method. Our findings are displayed in Tab.~\ref{tab:1}.
 
 Interestingly, all eigenvalues appear to be real, even though there is no a priori reason for that. Actually, following the reasoning in \cite{BiernatBizonMaliborski2016}, Sec.~3, one can prove that the eigenvalues for which $\text{Re}\:\lambda >(d-3)/2$  are necessarily real. Note, however, that this is of little significance as already for $d=9$ we observe no eigenvalues in this region. Note also that in each dimension, in addition to $\lambda_0=1$ there is exactly one more non-negative eigenvalue. However, it is very difficult to rigorously prove this observation. Nonetheless, for $d=7$ this eigenvalue happens to be $\lambda_1=3$ with an explicit corresponding eigenfunction, and we were in fact able to prove the following result.

 	\begin{proposition}[\cite{GlogicSchoerkhuber}, case $\ell=0$]
 		Let $d=7$. For $\mathrm{Re\,} \lambda \geq 0$,  the only solutions $f \in C^{\infty}[0,1]$ to Eq.~\eqref{Eq:ODE_L0} are (constant multiples of)
 		\begin{equation}
 		f_0(\rho) =  \frac{1-\rho^{2}}
 		{(1+\rho^{2})^{2}} \quad \text{and} \quad f_1(\rho) = \frac{1}{(1+\rho^2)^2}
 		\end{equation}
 		for which $\lambda=1$ and $\lambda=3$ respectively.
 	\end{proposition}

 \begin{table}[t]
 	\setlength{\tabcolsep}{10pt}
 	\centering
 	\begin{tabular}{clllll}
 		\hline\hline
 		$d$ & $\lambda_{1}$ & $\lambda_{0}$ & $\lambda_{-1}$ & $\lambda_{-2}$ & $\lambda_{-3}$ \\ \hline
 		5  & 4.37213 & 1.00000 & -0.53721 & -1.88858 & -3.17611 \\
 		6  & 3.39524 & 1.00000 & -0.54896 & -1.96235 & -3.32602  \\
 		7  & 3.00000 & 1.00000 & -0.55242 & -2.00000 & -3.41077  \\
 		8  & 2.78200 & 1.00000 & -0.55388 & -2.02356 & -3.46699  \\
 		9  & 2.64296 & 1.00000 & -0.55462 & -2.03986 & -3.50744 \\
 		\vdots & \vdots & \vdots & \vdots &  \vdots & \vdots \\
 		$\infty$ & 2.00000 & 1.00000 & -0.55593 & -2.13344 & -3.76974 \\
 		\hline\hline
 		\vspace{1ex}
 	\end{tabular}
 	\caption{All eigenvalues of Eq.~\eqref{Heun} appear to be real. In addition to two non-negative eigenvalues, there seem to be infinitely many negative ones, the largest three of them being shown.}
 	\label{tab:1}
 \end{table}
 
 On the other hand, it seems that there are infinitely many negative eigenvalues; in Tab.~\ref{tab:1} we listed the largest three of them. Here we should point out that (for $d=7$) only the eigenvalues for which $\text{Re\,}\lambda > -1/2$ correspond to isolated eigenvalues of the operator $\mathbf{L}_0$, see Lemma 5.2 in \cite{GlogicSchoerkhuber}. 
 This is related to the fact that we study Eq.~\eqref{Evol_perturb} in $H^3 (\mathbb{B}^7) \times H^2 (\mathbb{B}^7)$ and the spectral cut-off $-1/2$ is dictated by the choice of (the regularity of) this space. What is more, by increasing the Sobolev exponents the spectral cut-off is pushed to the left and in this way more and more negative eigenvalues get uncovered.

 A particularly interesting feature of our numerical results in Tab.~\ref{tab:1} is that eigenvalues seem to decrease as the dimension increases. It is therefore natural to ask as to whether they have limiting values as $d$ goes to infinity. To further investigate this, we rescale the independent variable
 \begin{equation}\label{scaling}
 x=dz
 \end{equation}
 in Eq.~\eqref{Heun} and by letting $d$ go to infinity we obtain the following equation
 \begin{align}\label{limiting_equ}
 \begin{split}
 zy''(z)+\left(-\frac{1}{2z}+\frac{4}{3z+1}+\lambda-\frac{3}{2} \right)y'(z)+\frac{\lambda-2}{4}\frac{3(\lambda-3)z + \lambda+5}{z(3z+1)}y(z)=0.
 \end{split}
 \end{align}
 In this process the interval $[0,1]$ contracts into a single point $z=0$. We therefore look for solutions of Eq.~\eqref{limiting_equ} that are analytic at $z=0$ only. There is indeed a formal power series solution centered at zero, which is however generically not convergent as zero is an irregular singular point. Therefore, by requiring convergence we impose a quantization condition on $\lambda$. This is however a different problem from the one above (where the relevant solutions are a priori analytic) and we therefore need a different approach. For this, we use an adaptation of the so-called continued fraction method.  This method relies on a remarkable observation which was (to the best of knowledge of the authors) first time made by Jaff\'e in \cite{Jaf34}. He makes use of the connection between the three term recurrence relations and continued fractions to compute the bound states of the hydrogen atom. This technique was later popularized inside the general relativity community by Leaver \cite{Lea85}, who used it to calculate the quasinormal modes of Kerr black holes. Interestingly, as pointed out by Bizo\'n \cite{Biz05}, although relatively old, this method does not seem to be well-known within the mathematics community.

 We start with the normalized formal power series solution to Eq.~\eqref{limiting_equ} centered at zero
 \begin{equation}\label{formal_series}
 y(z)=\sum_{n=0}^{\infty}a_n(\lambda) z^n.
 \end{equation}
 Here the coefficients are given by the recurrence relation
 \begin{equation}\label{Eq:rec_rel}
 a_{n+2}(\lambda)=A_n(\lambda) a_{n+1}(\lambda)+B_n(\lambda) a_n(\lambda)
 \end{equation}
 where
 \begin{align*}
 A_n(\lambda)=\frac{\lambda^2+(4n+7)\lambda+4n^2+8n-6}{2(n+2)}
 \end{align*}
 and
 \begin{align*}
 B_n(\lambda)=\frac{3(\lambda+2n-2)(\lambda+2n-3)}{2(n+2)}
 \end{align*}
 and the initial condition is $a_0(\lambda)=1$, $a_1(\lambda)=A_{-1}(\lambda)=(\lambda^2+3\lambda-10)/2$. The radius of convergence of the series \eqref{formal_series} is determined by the asymptotics of the coefficients $a_n$, and since this sequence satisfies Eq.~\eqref{Eq:rec_rel} we make use of the difference equation theory to determine all possible scenarios. For this we refer the reader to an excellent book by Elaydi \cite{Ela05}. First, it can be proved that for every $\lambda\in \mathbb{C}$ there are two linearly independent solutions to Eq.~\eqref{Eq:rec_rel} with the following behavior
 \begin{equation}\label{Eq:asymp}
 a_n^{(1)}(\lambda) \sim n! \, 2^n \,  n^{\lambda-\frac{1}{2}} \quad  \text{and} \quad a_n^{(2)}(\lambda) \sim (-3)^n \, n^{-6}
 \end{equation}
 as $n\rightarrow \infty$. This, together with the fact that $\lim_{n\rightarrow \infty}A_n(\lambda)=\lim_{n\rightarrow}B_n(\lambda)=\infty$, implies that for every $\lambda\in\mathbb{C}$ there are $c_1(\lambda)$ and $c_2(\lambda)$ such that
 \begin{equation}\label{Eq:comb}
 a_n(\lambda)=c_1(\lambda)a_n^{(1)}(\lambda) + c_2(\lambda)a_n^{(2)}(\lambda) 
 \end{equation}
 for large $n$. Also, note that since $a_n^{(1)}$ is the dominant solution the choice of the coefficient $c_1(\lambda)$ is unique, i.e.~it does not depend on the particular choice of the solution $a^{(1)}_n$.

 Now from Eqs.~\eqref{Eq:comb} and \eqref{Eq:asymp} it is clear that the eigenvalues are precisely the zeros of the function $c_1$. In the language of difference equations theory this is equivalent to saying that $a_n(\lambda)$ is a minimal solution to Eq.~\eqref{Eq:rec_rel}. Here minimal solution denotes the one that is asymptotically negligible relative to any other solution that is not its constant multiple (for large values of $n$). 
 Now, consider the following continued fraction
 \begin{equation}\label{Eq:continued_frac}
 \frac{B_0(\lambda)}{A_0(\lambda)+}\frac{B_1(\lambda)}{A_1(\lambda)+}\frac{B_2(\lambda)}{A_2(\lambda)+}\dots,
 \end{equation}
 written in the standard compact form. Pincherle's theorem, see \cite{Ela05}, p.~402, says that \eqref{Eq:continued_frac} is convergent if Eq.~\eqref{Eq:rec_rel} has no minimal solution with initial condition $a_0=0$, $a_1=1$. In this case, by Pincherle's theorem again, a solution to Eq.~\eqref{Eq:rec_rel} is minimal if and only if the ratio $a_1/a_0$ equals \eqref{Eq:continued_frac}. Therefore, the zeros of the following transcendental function
 \begin{equation}\label{Eq:transc_funct}
 A_{-1}(\lambda)+\frac{B_0(\lambda)}{A_0(\lambda)+}\frac{B_1(\lambda)}{A_1(\lambda)+}\frac{B_2(\lambda)}{A_2(\lambda)+}\dots
 \end{equation}  
 are eigenvalues. Furthermore, the roots of this function can be approximately computed by truncating the continued fraction for some large $n$ and then solving the corresponding equation. We displayed our results in the last row of Tab.~\ref{tab:1}. 
 
The preceding conclusion was made under the assumption that $\lambda$ is not such that the solution to Eq.~\eqref{Eq:rec_rel} with $a_0=0$, $a_1=1$ is minimal. As a matter of fact, such values of $\lambda$ exist and they can be computed by truncating \eqref{Eq:continued_frac} for large $n$ and then calculating the zeros of the reciprocal of the resulting function. What is more, it appears that there are infinitely many of these values and in general it is possible that one of them coincides with an eigenvalue causing its ``loss", in a sense that it then fails to be a zero of the function \eqref{Eq:transc_funct}. However, this does not seem to happen in our case.

Interestingly, Tab.~\ref{tab:1} suggests that one eigenvalue of Eq.~\eqref{limiting_equ} is $\lambda_1=2$. In fact, this is obvious since any constant function is a solution in that case. Furthermore, $\lambda_0=1$ is an eigenvalue with eigenfunction $y_0(z)=1-3z$, as can be found by calculating the limit of \eqref{gauge_mode}  as $d \rightarrow \infty$, having in mind \eqref{Change_var_Heun} and \eqref{scaling}. Then, since both (numerically observed) non-negative eigenvalues of Eq.~\eqref{limiting_equ} are known explicitly together with their corresponding eigenfunctions it is likely that an adjustment of spectral techniques developed in \cite{CosDonGloHua16}, \cite{CosDonGlo17} and \cite{GlogicSchoerkhuber} would yield a rigorous proof of non-existence of other unstable eigenvalues. Subsequently, by using some kind of perturbative argument to prove that for large values of $d$ the spectrum of Eq.~\eqref{Heun} is close to the one of Eq.~\eqref{limiting_equ} one would obtain an analog of Theorem \ref{Codim_one_radial} for all large dimensions $d$.

We note that the continued fraction method can also be used to compute the spectrum of Eq.~\eqref{Heun}. We checked that it yields the same result as displayed in Tab.~\ref{tab:1}.

 \smallskip
 
Finally, we remark that in the case $d=7$ the genuine instability $\lambda_1=3$ can be related to the conformal invariance of Eq.~\eqref{Lin_Eq_phys} which is due to the self-similar character of the potential. In fact, it is easy to check that the transformation $\phi \mapsto \tilde \phi$,
\begin{equation}\label{ConfTrans_T}
	\tilde \phi(t,r)  = ((T-t)^2 -r^2)^{-\frac{d-1}{2}} \phi \left (\tfrac{t-T}{(T-t)^2 - r^2} + T, 	\tfrac{r}{(T-t)^2 - r^2} \right ),
\end{equation}
leaves Eq.~\eqref{Lin_Eq_phys} invariant. In particular,
\begin{equation}
	\tilde \phi_0(t,r) = (T-t)^{3 - d} (1 - \tfrac{r}{T-t} )^{\frac{5-d}{2}} f_0(\tfrac{r}{T-t})
\end{equation}
gives rise to a solution of Eq.~\eqref{Evol_perturb_lin} given by
\begin{equation}
	\tilde \varphi_{0}(\tau,\rho) = e^{(d-4)\tau} (1 - \rho^2)^{\frac{5-d}{2}} f_0(\rho).
\end{equation}
For $d=5$, this is just an identity. In higher space dimensions, this transformation induces a singularity at the boundary of the cylinder, since  $f_0(\rho) = \mathcal O(1)$ around $\rho =1$ in general.  However, in $d=7$, $f_0(1) =0$ and thus
\begin{equation}
	\tilde \varphi_{0}(\tau,\rho) = e^{3 \tau} (1 - \rho^2)^{-1} f_0(\rho)
\end{equation}
corresponds to the mode solution for $\lambda_1 = 3$. This effect can be observed also in other supercritical wave equations, see e.g.~\cite{BiernatBizonMaliborski2016} for the wave maps problem.

\section{Threshold behavior}
\label{sec:Non-genericBlowup}

To understand the nature of the threshold for blowup we study the time-evolution from a numerical point of view. The simulations were performed by the second author.  We use two different numerical approaches.
The first one is very efficient for numerical studies of self-similar blowup and follows the scheme introduced in \cite{BiernatBizonMaliborski2016}.
It uses a suitable redefinition of variables and provides very accurate results for solutions which exhibit self-similar blowup.
Although being very efficient for observing blowup, it is not suitable to simulate solutions that disperse.
Therefore we solve the equation using the original $(t,r)$ variables to provide evidence for dispersion of subcritical initial data.

In the following, we present results for $d=5$ and $d=7$ as an illustration. However, we expect analogous behavior for higher dimensions (we did simulations in $d=9$ and $d=11$ indicating this).

\subsection{Numerical approach}
\label{sec:NumericalApproach}

\begin{figure}[t]
  \centering
  \includegraphics[width=0.475\textwidth]{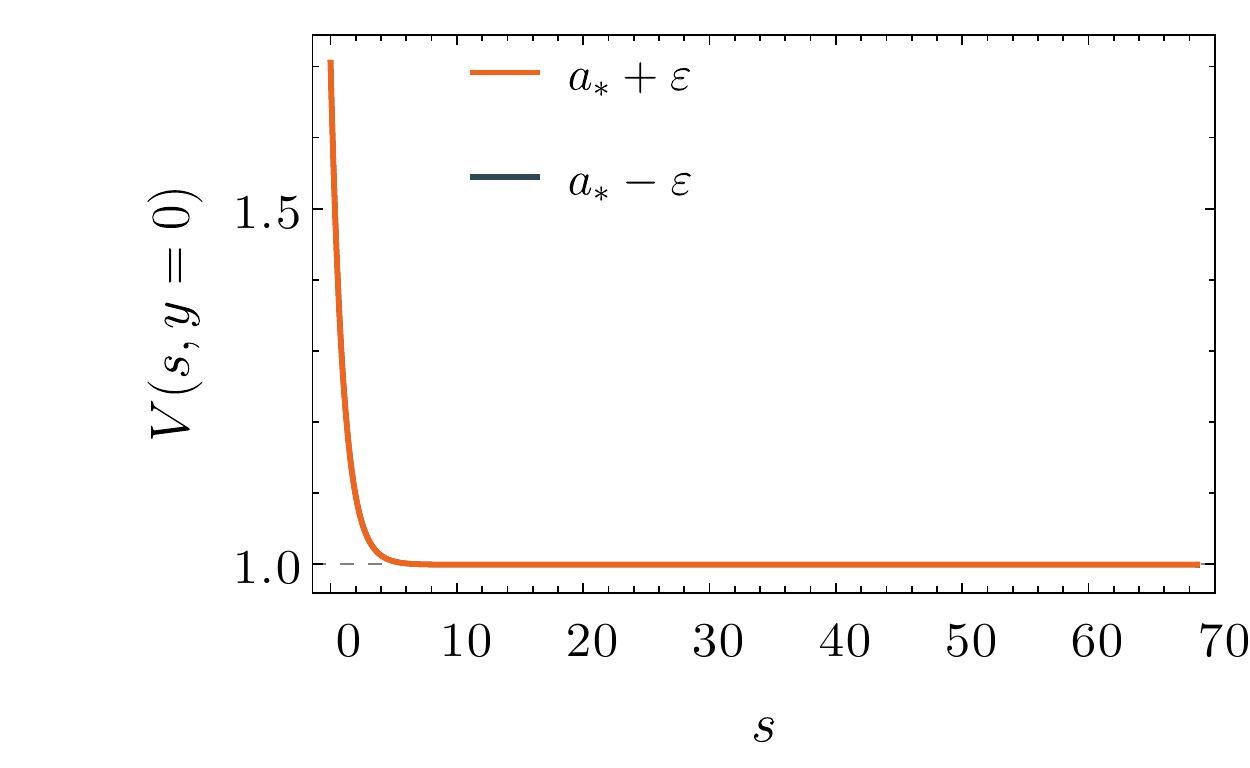}
  \includegraphics[width=0.475\textwidth]{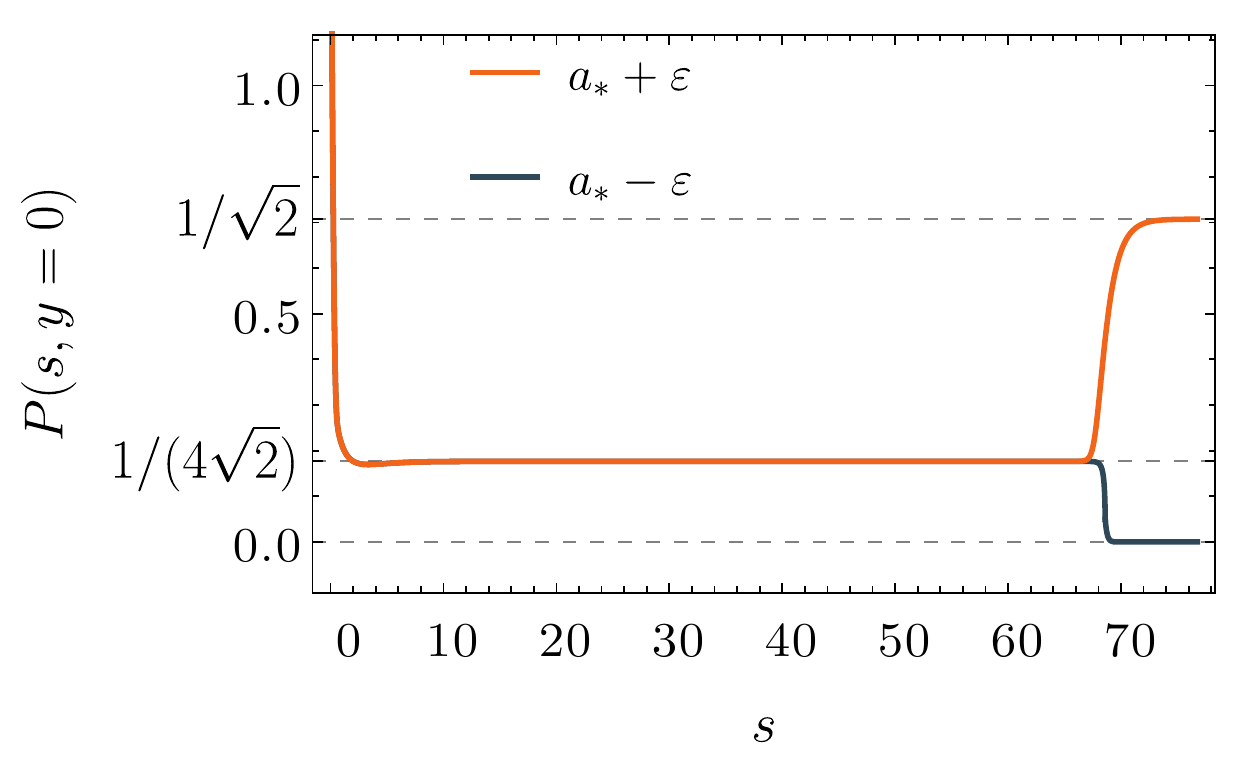}
  \caption{The evolution of marginally sub- (blue line) and
    supercritical (orange line) initial data given by Eq.~\eqref{eq:16},  in $d=5$ in computational
    variables (\ref{eq:11})-(\ref{eq:12}).}
  \label{fig:1}
\end{figure}

Following \cite{BiernatBizonMaliborski2016} we introduce new (computational) coordinates $(s,y)$
through
\begin{equation}
  \label{eq:11}
  t=\int e^{-s}h(s)\diff{s}, \quad
  r=e^{-s}y.
\end{equation}
The function $h$ introduced in \eqref{eq:11} is used to make the coordinate
transformation adapt to the blowing up solution; note that $(s,y)$ are
the self-similar coordinates $(\tau,\rho)$ defined in
(\ref{SS_coordinates}) when $h(s)\equiv 1$. We also define new
dependent variables
\begin{equation}
  \label{eq:12}
  V(s,y) = e^{-s}u(t,r), \quad
  P(s,y) = e^{-2s}\partial_{t}u(t,r).
\end{equation}
The time evolution of the new variables (\ref{eq:12}) follows from
(\ref{eq:11}) and the equation (\ref{eq:1}); explicitly we have
\begin{eqnarray}
  \label{eq:13}
  \begin{split}
    \partial_{s}V(s,y) &= h(s)P(s,y) - V(s,y) - y\partial_{y}V(s,y), \\
    \partial_{s}P(s,y) &= h(s)\left(\partial_{y}^{2}V(s,y) +
      \frac{d-1}{y}\partial_{y}V(s,y) + V(s,y)^{3}\right) \\
    & \ \ \ - 2P(s,y) - y\partial_{y}P(s,y).
  \end{split}
\end{eqnarray}
The main advantage of this rescaling and the redefinition of dependent
variables is that $V$ stays finite if we set
$h(s)=1/P(s,0)$. More precisely, this choice of $h$ leads to the following long time asymptotics at $y=0$: $V(s,0) = 1 + C e^{-s}$,
for some constant $C\in\mathbb{R}$. Most importantly, $P(s,0)\rightarrow 1/U(0)$ in the case of self-similar blowup, where $U$ is a self-similar profile. Note that this behavior, finiteness of the dependent variables, is in stark contrast to the original variables
$(u,\partial_t u)$ which blow up in finite time.  Some care has to be taken to assure that \eqref{eq:11} defines a coordinate transformation, in particular for our choice of $h$ it is evident that $P(s,0)=0$ for finite $s$ is problematic. It particular, this happens when $\partial_t u(t ,0)$ vanishes for some $t$. This is the reason why these coordinates are not suitable for the study of solutions that exhibit oscillating behavior. We overcome this by choosing initial conditions which avoid this behavior, as it is apparent from Fig.~\ref{fig:1}

We also note that in order to speed up the numerical calculations we choose the initial data for which the decision on either blowup or dispersion can be made relatively early so that we optimise for computational resources. Other than that the initial data we use are generic. We note that the evolution is stopped when $P(s,0)$ reaches a small positive value. To understand the behavior beyond this point, we evolve the corresponding initial data in physical coordinates $(t,r)$ and observe dispersion.


As pointed out in \cite{BiernatBizonMaliborski2016} the transformation
(\ref{eq:11}) with a proper choice of $h$ introduces self-adapting
coordinates that accurately resolve both spatial and temporal scales
of blowing up solution. In this way we avoid using an adaptive spatial
mesh and a time rescaling techniques. In fact to solve (\ref{eq:13})
we use a standard method of lines with 6th order finite difference
approximation in space and a 6th order Runge-Kutta method with fixed
time step as a time stepping algorithm. Staggered spatial grid with
fixed mesh size is used to deal with the $y=0$ singularity in
(\ref{eq:13}). In addition, a symmetry of the variables (\ref{eq:12})
at the origin is used to construct finite difference stencils close to
the coordinate singularity. We add a standard dissipation term to
suppress high frequency noise in the data introduced by spatial
discretisation. Together with (\ref{eq:13}) we solve the differential of (\ref{eq:11}) with the initial condition $t(0)=0$ to get the relation $t(s)$ needed for the subsequent analysis. The code was written in Mathematica whose flexibility
and functionality allowed us to use arbitrary precision arithmetics
seamlessly. We remark that higher precision is crucial to get close enough to
the critical solution and obtain a detailed description of near
critical evolutions. To speed up computations we parallelized the
bisection search (discussed in the following section) probing the
search interval using multiple (typically 64) cores simultaneously.

\subsection{Results}\label{sec:Results}

For the family of initial data
\begin{equation}
  \label{eq:16}
  V(0,y)=P(0,y)=\frac{a}{\cosh{y}},
\end{equation}
we integrate \eqref{eq:13} forward in $s$ with $a > 0$ as the only free parameter. For large amplitudes $a$ the solution approaches the ODE blowup \eqref{Eq:ODEblowup}. In particular, $P(s,0)$ goes to $1/U_{0}(0)=1/\sqrt{2}$ as $s\rightarrow\infty$, see Fig.~\ref{fig:1}. We perform a bisection search in the amplitude $a$ based on the criteria outlined in the previous section; in this way, we find the critical value of $a$, which we call $a_*$, at threshold. 
For this value, $P(s,0)$ is approximately $1/U^*(0)$, which indicates that $U^*$ is an attractor within the threshold. Explicitly, $a_{*}\approx 1.710572581$ in $d=5$
and $a_{*}\approx 2.335609125$ in $d=7$.

\begin{figure}[t]
  \centering
  \includegraphics[width=0.75\textwidth]{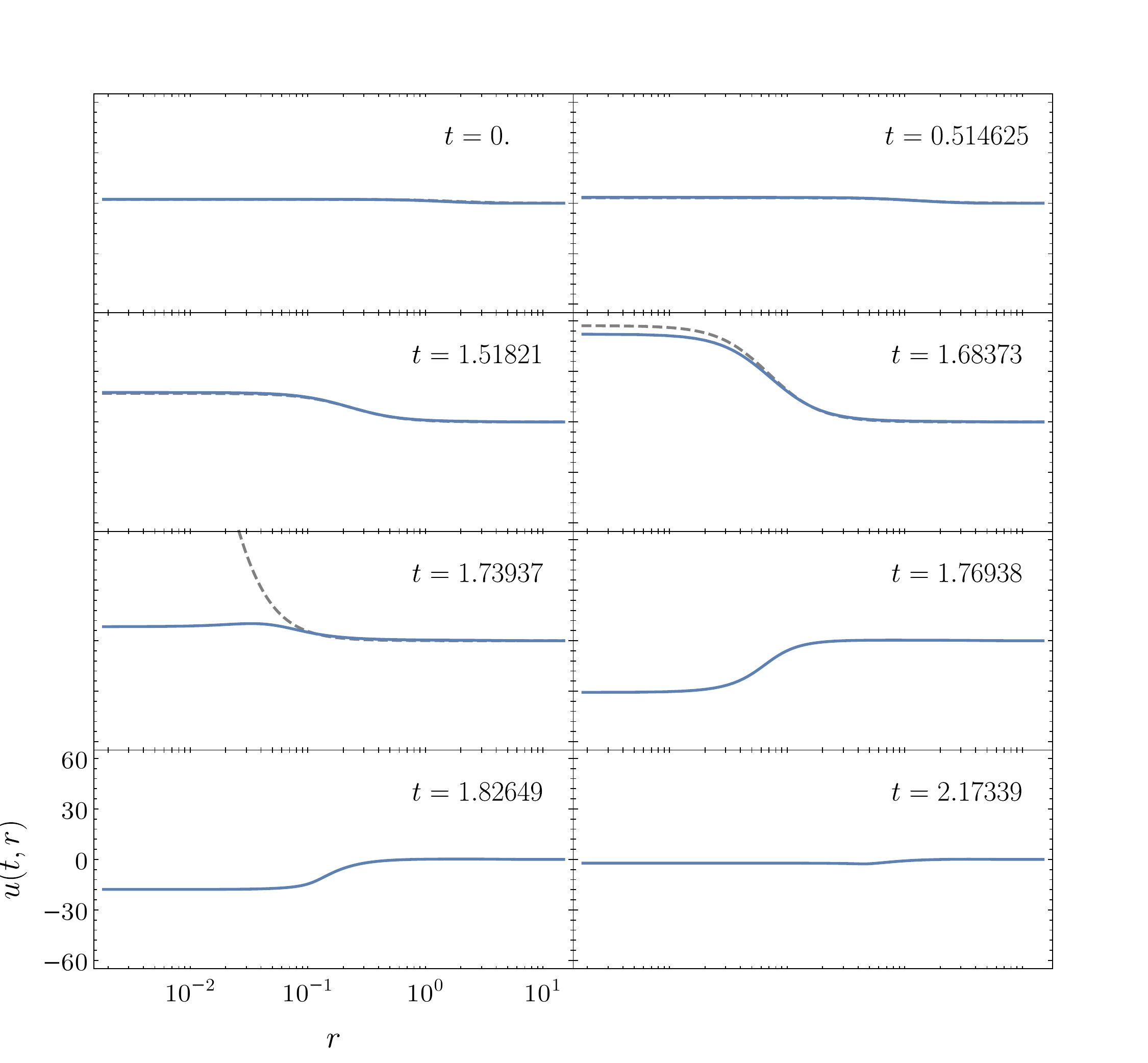}
   \caption{The evolution of subcritical initial data for $d=7$ in $(t,r)$. Blue line shows the numerical solution with the initial amplitude $a=2.33558<a_{*}$, which for intermediate times approaches $u_T^*$ (dashed line) with the blowup time fixed to $T\approx 1.7536$.}
  \label{fig:2}
\end{figure}

\begin{figure}[t]
  \centering
  \includegraphics[width=0.37\textwidth]{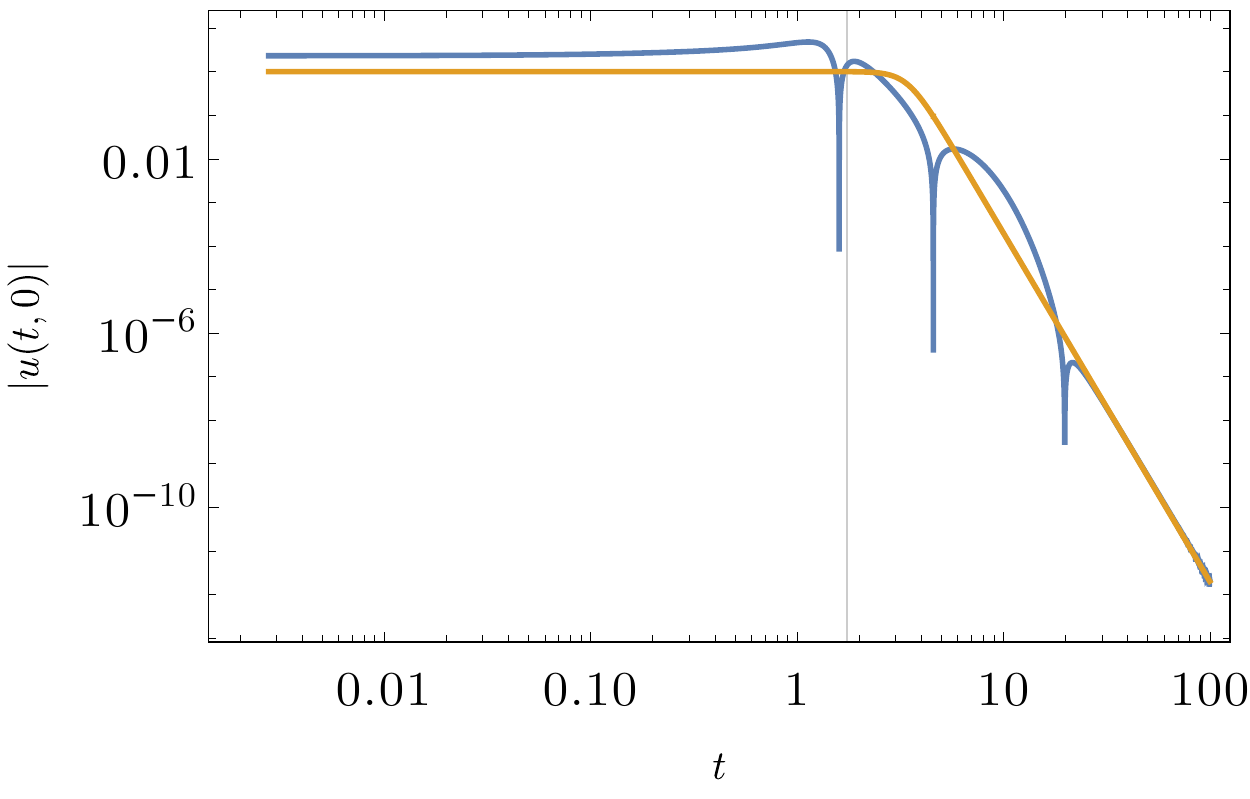}
  \hspace{4ex}
  \includegraphics[width=0.45\textwidth]{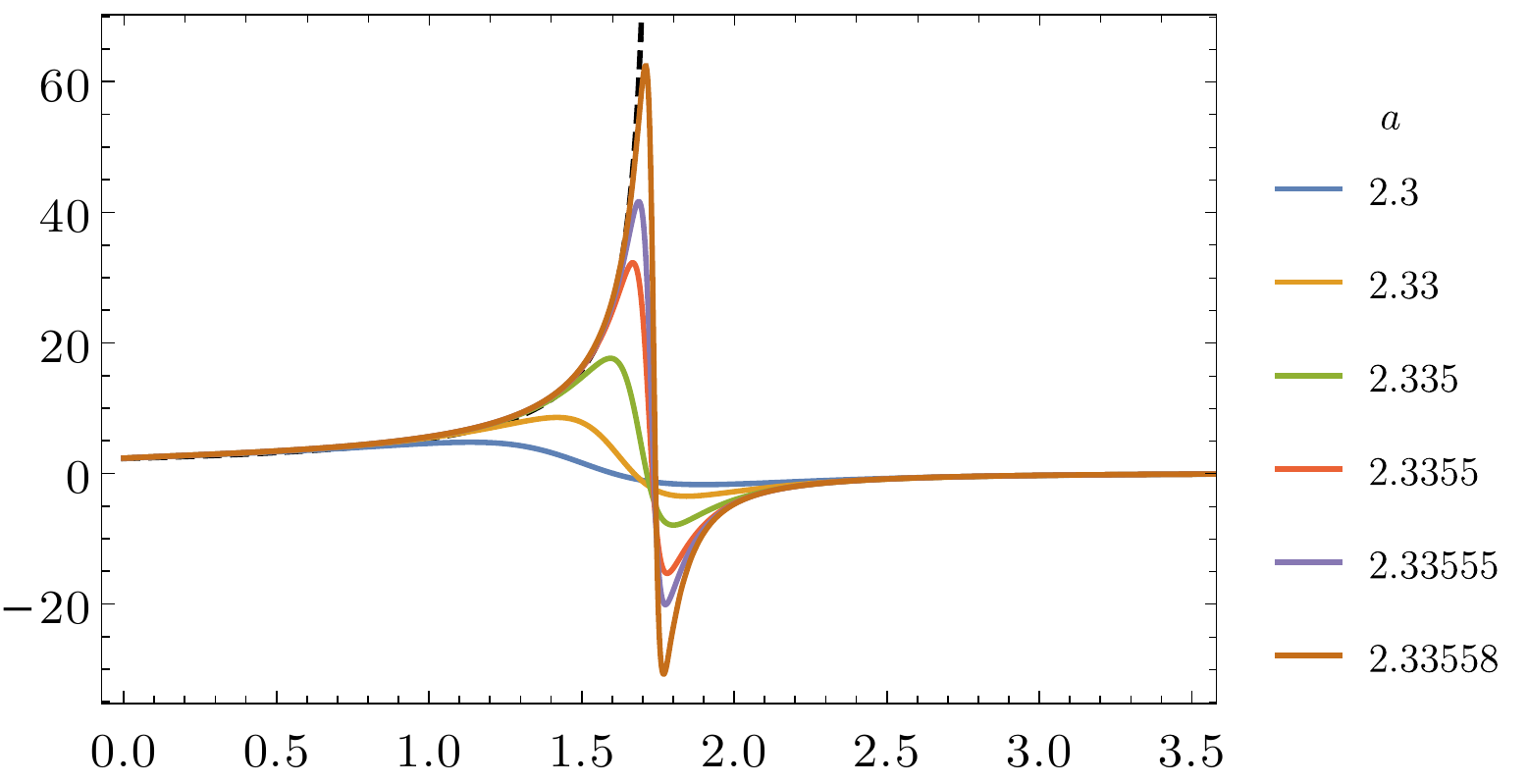}
  \caption{ 
  For $d=7$ subcritical  initial data ($a = 2.3 <a_{*}$) disperse via a polynomial decay $t^{-8}$  (left panel).
 When $a \rightarrow (a_{*})^{-}$ the solution approaches the self-similar solution $u_T^*$ (right plot, dashed line). The closer the amplitude (color coded) to the critical amplitude the longer the solution follows $u_T^*$ . Ultimately, all subcritical solutions decay to zero.}
  \label{fig:3}
\end{figure}

To analyze the subcritical evolution, we evolve initial data for $a < a_*$ in $(t,r)$. We solve (\ref{eq:1}) using the same algorithms as used for solving the system (\ref{eq:13}). We note that at the initial time $s=t=0$ the transformation in (\ref{eq:11})-(\ref{eq:12}) implies that $r=y$, and 
\begin{equation}
	\label{eq:22XII19_01}
	V(0,y) = u(0,r), \quad P(0,y)=\partial_{t}u(0,r)\,.
\end{equation}
To see the intermediate attractor at the threshold we have to find the approximate blowup time of the critical evolution corresponding to $a \approx a_*$. This procedure is given in \cite{BiernatBizonMaliborski2016}. By this, we observe that subcritical data approaches $u_T^*$ for intermediate times before it disperses, see Fig.~\ref{fig:2}.
Furthermore, the behavior at the origin for different values of $a < a_*$ is displayed in Fig.~\ref{fig:3}. 

\begin{figure}[t]
  \centering
  \includegraphics[width=0.75\textwidth]{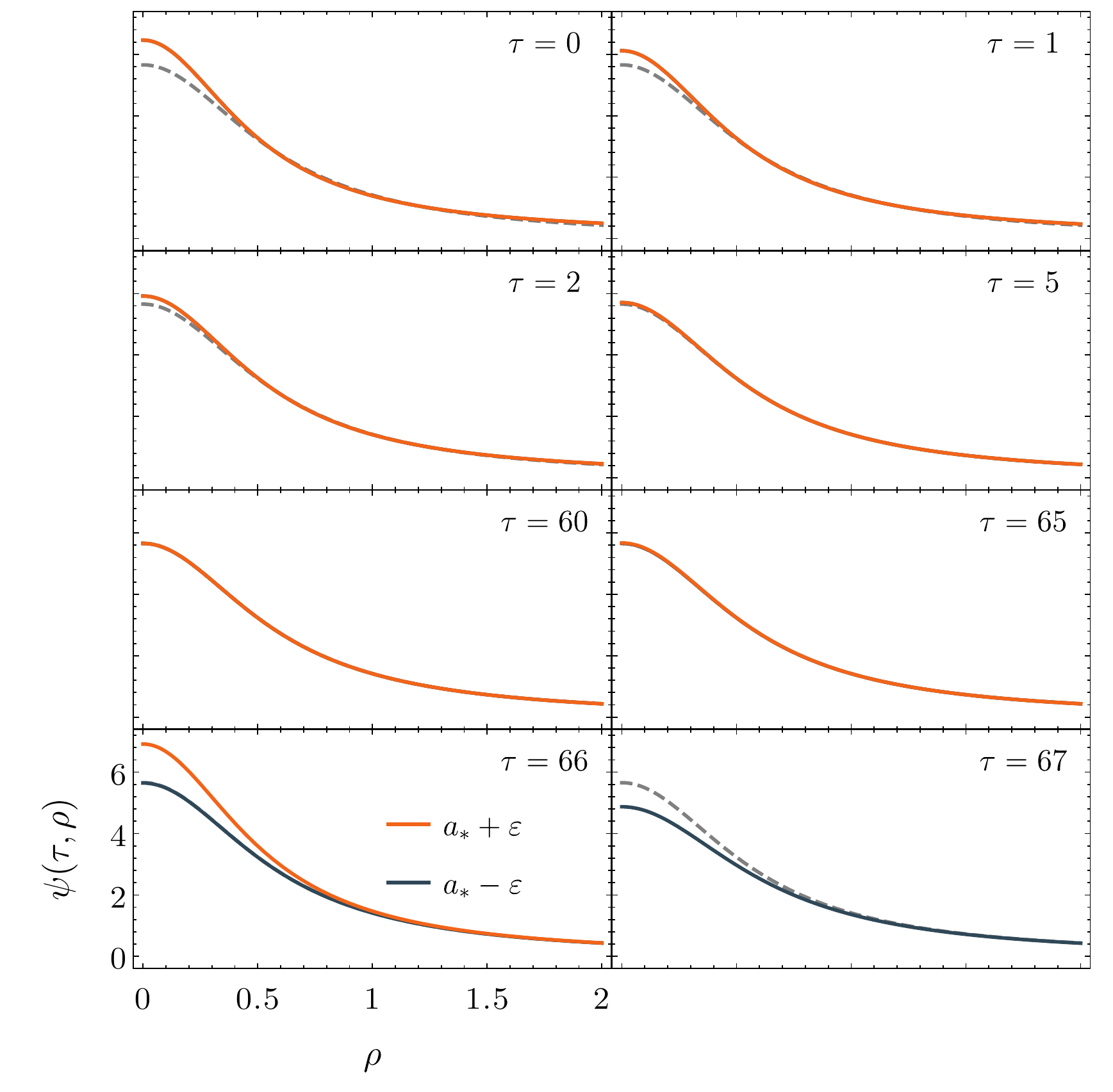}
  \caption{The evolution of marginally sub- (blue line) and
    supercritical (orange line) evolutions for $d=5$ in self-similar
    coordinates adapted to the threshold solution. Both solutions
    approach the intermediate attractor $U^{*}$ (dashed line). 
    After some time solutions depart from $U^{*}$ in opposite
    directions. In the last frame, the supercritical solution is out
    of range of the plot as the coordinates used differ from the ones
    in which we would see the approach to $U_{0}$ (the generic
    blowup).}
  \label{fig:4}
\end{figure}

\begin{figure}[t]
  \centering
  \includegraphics[width=0.75\textwidth]{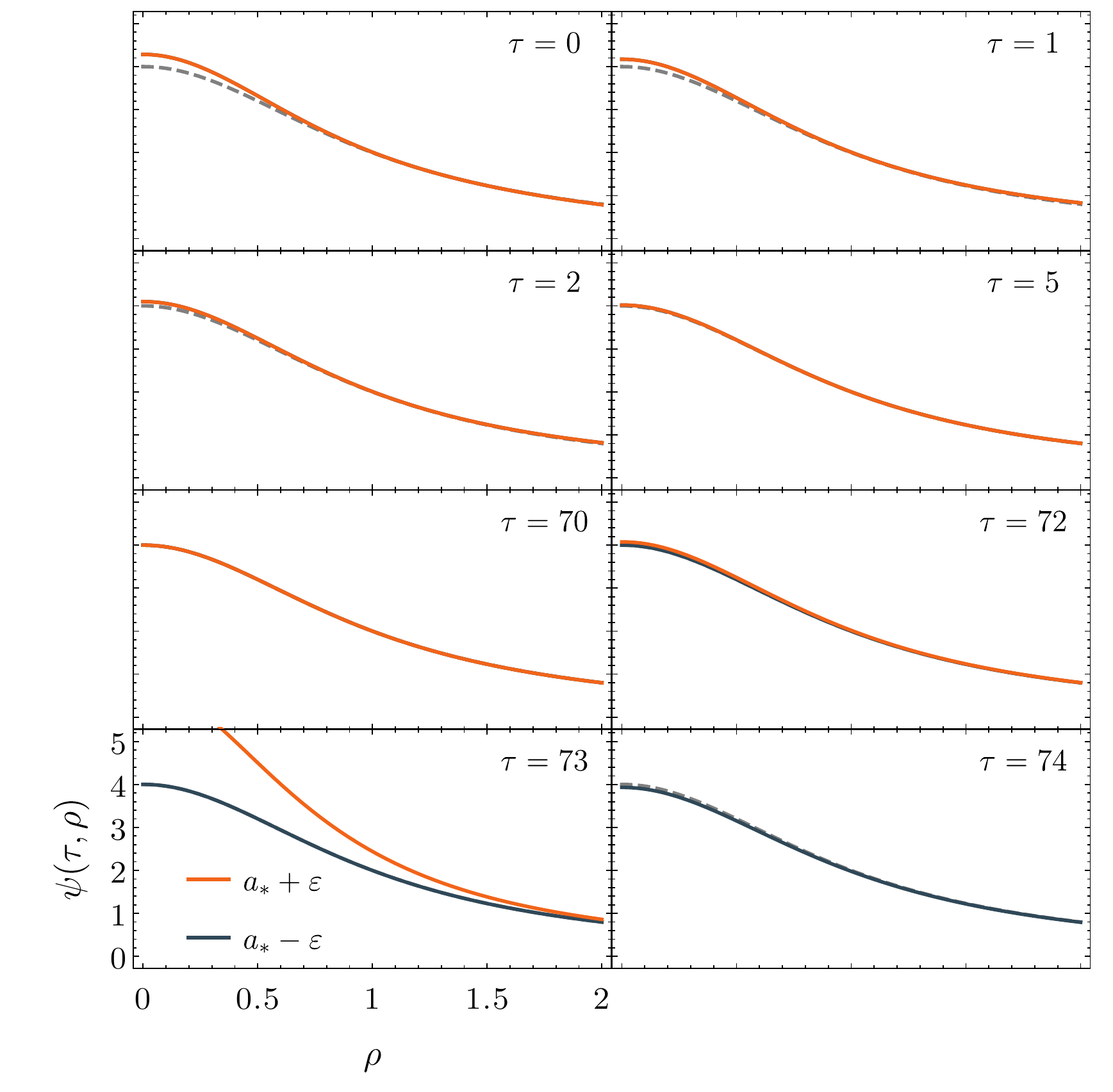}
  \caption{Same as Fig.~\ref{fig:4} but for $d=7$.}
  \label{fig:5}
\end{figure}

As an additional evidence for the threshold nature of $u_T^*$, we pass to self-similar coordinates $(\tau, \rho)$ using the computed critical blowup time, see Fig.~ \ref{fig:4} and Fig.~ \ref{fig:5}.

\medskip
Finally, we perform a qualitative comparison of the numerical data with the analytical predictions of Sec.~\ref{Sec:ODE_Analysis}. We expect a solution in self-similar coordinates, recall Eq.~\eqref{SS_coordinates2}, to behave like
\begin{equation}
  \label{eq:17}
  \psi(\tau, 0) = c + a_{1}e^{\lambda_{1}\tau} + a_{0}e^{\tau} +
  a_{-1}e^{\lambda_{-1}\tau} + \cdots,
\end{equation}
close to the critical solution, where the dots stand for
even faster decaying modes. The bisection procedure ensures that the
coefficient of the unstable mode $a_{1}\sim a-a_{*}\equiv\varepsilon$
is small\footnote{Using 128 digits of precision for $d=5$ we were able
  to fine-tune $a$ up to $\varepsilon\approx 10^{-129}$ whereas for
  $d=7$ and higher with 96 digits we obtain $\varepsilon\approx 10^{-96}$.} so
that for long enough time we see the convergence of $\psi(\tau, 0)$ to $U^{*}(0)$.
We find that the amplitude of the gauge mode is typically larger than the amplitude of the genuine unstable mode. This is due to propagation of errors on different stages of data analysis. We stress that the appearance of the gauge mode is an artefact of the data preprocessing and is due to the uncertainty of the blowup time $T$ which is then used in the coordinate transformation.  

\begin{figure}[t]
  \centering
  \includegraphics[width=0.475\textwidth]{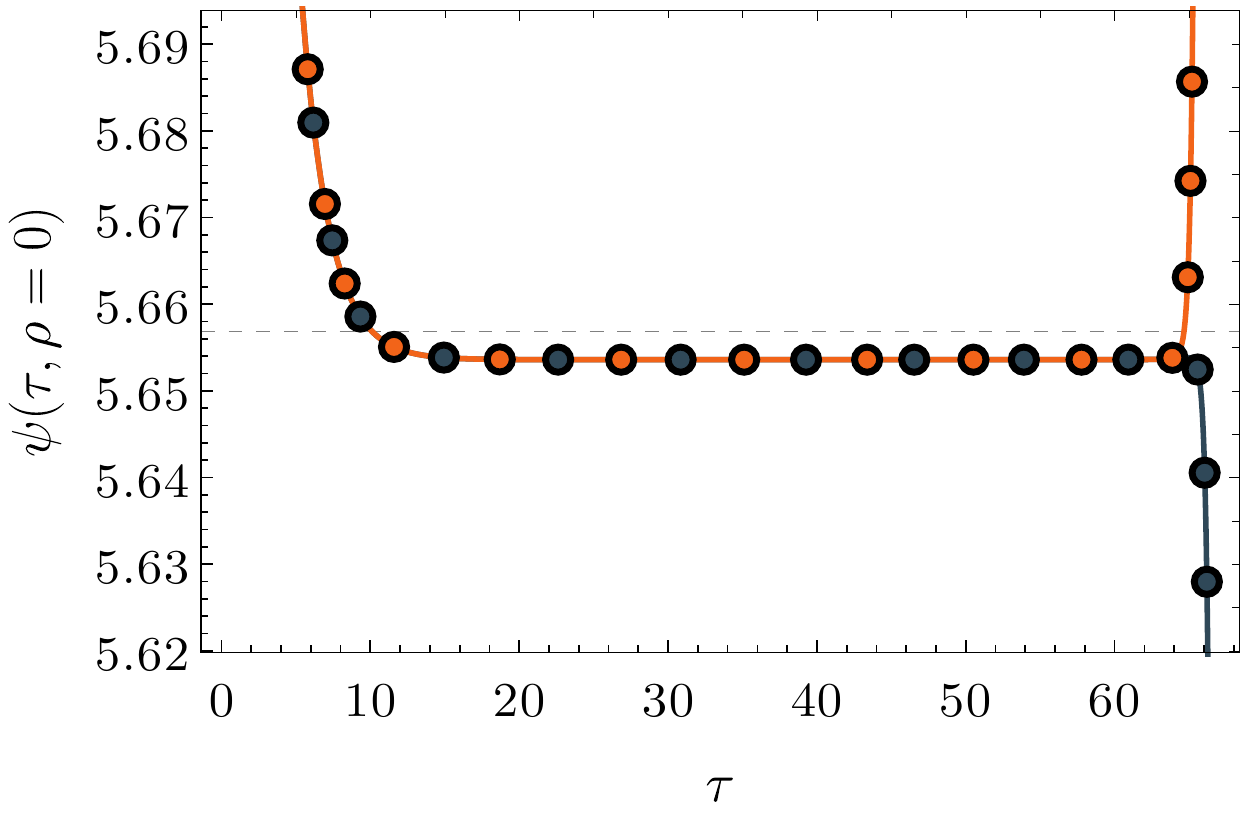}
  \hspace{2ex}
  \includegraphics[width=0.475\textwidth]{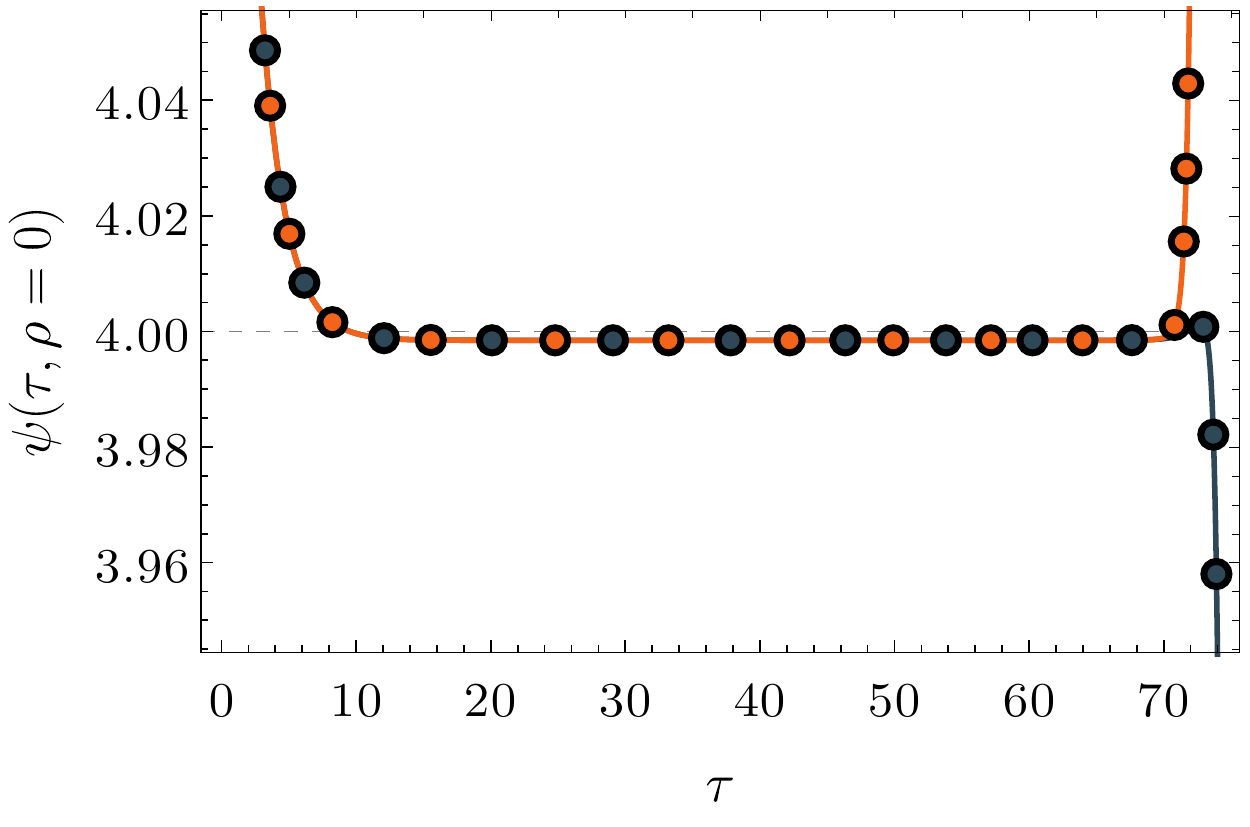}
  \caption{The convergence of nearly critical evolutions (blue line
    for subcritical case, orange line for supercritical case) to $U^{*}(0)$. The rate of convergence is determined by the
    stable eigenvalue $\lambda_{-1}$. At a later time, the unstable
    mode takes over and we see the divergence with the rate given by $\lambda_{1}$. 
    The gauge mode is removed by a proper choice of the blowup time $T$. Both $d=5$ (left panel) and $d=7$ (right panel)
    are presented.}
  \label{fig:6}
\end{figure}

The comparison of the theoretical prediction
(\ref{eq:17}) with the numerical results is presented in
Fig.~\ref{fig:6}, while the results of the fits are collected in
Tab.~\ref{tab:2}. Note that for $c$ we obtain approximately $U^{*}(0)$. The magnitude of $a_{0}$ is small and $a_{1}$ is small and of opposite sign in super- and subcritical evolutions. Furthermore, $\lambda_{-1}$ is in accordance with the numerically computed value in Tab.~\ref{tab:1}.

\begin{table}[t]
  \setlength{\tabcolsep}{9pt} \centering
  \begin{tabular}{crr}
    \hline\hline
    $d=5$ & subcritical & supercritical \\ \hline
    $c$      & 5.653595752 & 5.653595793 \\
    $a_{1}$  & $-5.423965718\cdot 10^{-128}$ & $5.002789061\cdot 10^{-126}$ \\
    $a_{0}$  & $2.184193927\cdot 10^{-32}$ & $2.144231964\cdot 10^{-32}$ \\
    $a_{-1}$ & 0.7585737069 & 0.7585801575 \\
    $\lambda_{-1}$ & $-0.5387977976$ & $-0.5387995060$ \\
    \hline\hline
    \\
    \hline\hline
    $d=7$ & subcritical & supercritical \\ \hline
    $c$ & 3.998463960 & 3.998463963 \\
    $a_{1}$  & $-2.707332535\cdot 10^{-98} $ & $1.138265563\cdot 10^{-95}$ \\
    $a_{0}$  & $1.107517572\cdot 10^{-34}$ & $1.106153203\cdot 10^{-34}$ \\
    $a_{-1}$ & 0.3001361956 & 0.3001362363\\
    $\lambda_{-1}$ & $-0.5535087559$ & $-0.5535088248$ \\
    \hline\hline
    \hspace{1ex}
  \end{tabular}
  \caption{Results of the fit of theoretical prediction to the sub-
    and supercritical evolutions ($d=5$ and $d=7$ cases),
    see~(\ref{eq:17}). The value of $\lambda_1$ is fixed according to the results of Sec.~\eqref{Sec:ODE_Analysis}. Listed values are results of the fit. }
  \label{tab:2}
\end{table}


\bigskip
\emph{Acknowledgement.} We would like to thank Piotr Bizo\'n  and Roland Donninger for helpful discussions. Furthermore, we thank the anonymous referee for his/her useful comments and suggestions.
\pagestyle{plain}
\bibliographystyle{unsrt}

\bibliography{references.bib}

\end{document}